\documentclass{ACmod}

\newcommand{\bbA}{\mathbb{A}}

\newcommand{\bbD}{\mathbb{D}}
\newcommand{\bbH}{\mathbb{H}}
\newcommand{\bbP}{\mathbb{P}}
\newcommand{\bbZ}{\mathbb{Z}}
\newcommand{\cF}{\mathcal{F}}
\newcommand{\ii}{{\mathtt i}}

\DeclareMathOperator{\SL}{SL}
\DeclareMathOperator{\Imag}{Im}
\newcommand{\ra}{\to}
\newcommand{\btau}{{\bar{\tau}}}
\newcommand{\del}{\partial}
\newcommand{\chkF}{\check{F\,}\!}
\newcommand{\chkX}{\check{X\,}\!}
\newcommand{\pt}{{\mathsf{pt}}}
\newcommand{\Kah}{{\mathsf{K\ddot{a}h}}}
\newcommand{\cx}{{\mathsf{cx}}}
\newcommand{\PD}{{\mathsf{PD}}}
\renewcommand{\phi}{\varphi}

\begin{document}

\title{Moonshine at Landau-Ginzburg points}

\author[C\u ald\u araru, He and Huang ]{% 
Andrei C\u ald\u araru, Yunfan He, Shengyuan Huang}

\address{Mathematics Department,
University of Wisconsin--Madison, 480 Lincoln Drive, Madison, WI
53706--1388, USA\newline\normalfont{\texttt{andreic@math.wisc.edu,
  yhe233@math.wisc.edu, shuang279@wisc.edu}}}

\begin{abstract} 
  {\sc Abstract:} We formulate a conjecture predicting unexpected
  relationships among the coefficients of the elliptic expansions of
  Klein's modular $j$-function around $j=0$ and $j=1728$.  Our
  conjecture is inspired by recent developments in mirror symmetry, in
  particular by work of Tu~\cite{Tu} computing categorical enumerative
  invariants of matrix factorization categories and by work of
  Li-Shen-Zhou~\cite{LSZ} computing FJRW invariants of elliptic
  curves.
\end{abstract}

\maketitle

{\bf 1. The conjecture}
\setcounter{section}{1}

\paragraph
The Monstrous Moonshine conjecture describes a surprising
relationship, discovered in the late 1970s, between the coefficients
of the Fourier expansion of Klein's $j$-function around the cusp
\[ j(\tau) = \frac{1}{q} + 744 + 196884q + 21393760q^2 + 864299970q^3
  + 20235856256q^4 + \cdots\]
and dimensions of irreducible representations of the Monster group.
Fourier expansions of other modular forms around the cusp are
critically important in number theory and algebraic geometry.  In
particular such expansions appear directly in computations of
Gromov-Witten invariants of elliptic curves~\cite{Dij}.

\paragraph
In this note we study the elliptic expansion of the $j$-function
around the hexagonal point $j = 0$ and the square point
$j=1728$, instead of around the cusp $j=\infty$.  At $j=0$ the
elliptic curve is the Fermat cubic, cut out in $\bbP^2$ by 
$x^3 + y^3+ z^3 = 0$, while at $j=1728$ it is given by $x^4+y^4+z^2 =
0$ in the weighted projective space $\bbP^2_{1,1,2}$.

From an enumerative geometry perspective the fact that we work around
the hexagonal and square points instead of around the cusp suggests
that we are working with Fan-Jarvis-Ruan-Witten (FJRW) invariants
instead of Gromov-Witten invariants.  See~(\ref{subsec:justification})
for details.

\paragraph
Let $\bbH$ and $\bbD$ denote the upper half plane and the unit disk in
the complex plane, respectively.  Fix $\tau_* = e^{\pi\ii/3}$ or
$\tau_* = \ii$ as the points\footnote{Any other point in the
  $\SL(2,\bbZ)$ orbit of $\tau_*$ works equally well, with only minor
  changes in the constants below.}  in $\bbH$ around which to carry
out the expansion.
\newpage

The uniformizing map $S$ around $\tau_*$ is the map
\begin{align*}
  S : \bbH \ra \bbD, \quad\quad & S(\tau) = \frac{\tau-\tau_*}{\tau-\btau_*},
\intertext{with inverse}
  S^{-1}: \bbD \ra \bbH, \quad\quad & S^{-1}(w) = \frac{\tau_*-\btau_* w}{1-w}.
\end{align*}
The {\em elliptic} expansion of $j$ around $\tau_*$ is simply the
Taylor expansion of $j \circ S^{-1}$ around $w=0$.  Its
coefficients are closely related ~\cite[Proposition 17]{Zag} to the
values of the higher modular derivatives $\del^nj(\tau_*)$,
\[ j\left (S^{-1}(w)\right) = 
  \sum_{n=0}^\infty \frac{(4\pi\Imag \tau_*)^n \del^nj(\tau_*)}{n!} w^n.  \]

\paragraph
\label{subsec:rescale}
The values of the higher modular derivatives of $j$ can be computed
term-by-term by a well-known recursive procedure.  The results are
rational multiples of products of powers of the Chowla-Selberg
period\footnote{The  exact value of $\Omega$ is unimportant, but in
  this case $\Omega =  1/\sqrt{6\pi}\left(\Gamma(1/3)/\Gamma(2/3)
  \right)^{3/2}$ for the hexagonal point and $\Omega =
  1/\sqrt{8\pi}\left (\Gamma(1/4)/\Gamma(3/4)\right)$ for the square point.}
$\Omega$ and of $\pi$.

Let $s(w) = 2\pi\Omega^2 \cdot S(w)$ denote the rescaling of $S$ by
the factor $2\pi\Omega^2$.  Then around $\tau_* = \exp(\pi\ii/3)$ we have 
\[ j\left( s^{-1}(w) \right )  =
  13824 w^3 - 39744w^6 + \frac{1920024}{35}w^9 -
  \frac{1736613}{35}w^{12} + \cdots,  \]
while around $\tau_* = \ii$ we have
\[  j\left( s^{-1}(w) \right )  = 1728 + 20736w^2 +105984w^4
  +\frac{1594112}{5}w^6 + \frac{3398656}{5}w^8+\cdots .\]

\paragraph
\label{subsec:gh}
In his study of categorical Saito theory of Fermat
cubics~\cite[Section 4]{Tu} Tu introduced the following two power
series with rational coefficients: 
\begin{align*}
  g(t)  & = \sum_{n=0}^\infty (-1)^n
  \frac{\left((3n-2)!!!\right)^3}{(3n)!} t^{3n}, \\
  h(t)  & = \sum_{n=0}^\infty (-1)^n
  \frac{\left((3n-1)!!!\right)^3}{(3n+1)!} t^{3n+1}.
\end{align*}
He argued that the ratio $h(t)/g(t)$ gives a flat coordinate on the
moduli space of versal deformations $x^3+y^3+z^3+3txyz = 0$ of the
Fermat cubic.

Similarly, for the elliptic quartic we introduce the two power series below
\begin{align*}
  g(t)  & = \sum_{n=0}^\infty
          \frac{\left ((4n-3)!!!!\right)^2}{(2n)!} t^{2n}, \\
  h(t)  & = \sum_{n=0}^\infty
          \frac{\left((4n-1)!!!!\right)^2}{(2n+1)!} t^{2n+1}.
\end{align*}
Even though the notation $g,h$ appears overloaded, it should be
evident from context which power series we refer to.

Our main result is the following conjecture.

\begin{Conjecture}
  \label{conj:main}
  (a) Around the hexagonal point the elliptic expansion of the
  $j$-function satisfies 
  \begin{align*}
    j\left( s^{-1} \left ( \frac{h(t)}{g(t)}\right ) \right ) & 
                     = 27t^3\left (\frac{8-t^3}{1+t^3} \right )^3 \\
    & = 13824t^3 - 46656t^6 +  99144t^9 - 171315t^{12} + 263169t^{15}
      -\cdots .
  \end{align*}

  (b) Around the square point the elliptic expansion of the
  $j$-function satisfies
  \begin{align*}
    j\left( s^{-1} \left ( \frac{h(t)}{g(t)}\right )
    \right ) &  = (192 + 256t)\left ( \frac{3+4t}{1-4t^2} \right )^2 \\
    & = 1728 + 20736t^2 + 147456t^4 + 851968t^6 + 4456448t^8 +
      \cdots .
  \end{align*}
\end{Conjecture}

\paragraph {\bf Notes.}
It is remarkable that the coefficients in the above power
series are all integers, despite $j(s^{-1}(w))$ only having rational
coefficients.  We were unable to find other modular forms with this
property.  We verified the validity of the conjectures up to $t^{24}$
in both cases, by computer calculations.

\paragraph{\bf Acknowledgments.}
We would like to thank Junwu Tu, Jie Zhou, Michael Martens, and Ken
Ono for helping out at various stages of the project.

This work was partially supported by the National Science
Foundation through grant number DMS-1811925.

\newpage

{\bf 2. Mirror symmetry origin of the conjecture}
\setcounter{section}{2}
\setcounter{paragraph}{0}

\paragraph
The original statement of mirror symmetry is formulated as the
equality of two power series associated to a pair $(X, \chkX)$ of
mirror symmetric families of Calabi-Yau varieties.  These two power
series are
\begin{enumerate}
\item[(a)] the generating series, in a formal variable
  $Q$, of the enumerative invariants of the family $X$ (the
  A-model potential); 
\item[(b)] the Taylor expansion of a Hodge-theoretic function (the
  period) on the moduli space of complex structures $M^\cx$ of the
  mirror family $\chkX$, with respect to a flat coordinate $q$ on this
  moduli space (the B-model potential).
\end{enumerate}
In order to compare the two power series, the variables $q$ and $Q$
are identified via an invertible map $\psi$ called the mirror map.

In physics the formal variable $Q$ is viewed as a flat coordinate on
the (ill-defined mathematically) complexified K\"ahler moduli space
$M^\Kah$, and the mirror map is interpreted as an isomorphism
\[ \psi : M^\cx \ra M^\Kah \]
between germs of $M^\cx$ and $M^\Kah$ around special points.
Traditionally these special points are the large volume and large
complex limit points, respectively.

\paragraph
The original mirror symmetry computation of~\cite{COGP} follows this
pattern.  It predicts a formula for the generating series of genus
zero Gromov-Witten invariants of the quintic $X$, by equating it to the
expansion of a period (solution of the Picard-Fuchs equation) for the
family of mirror quintics $\chkX$.  The equality of the two sides
allows one to calculate the genus zero Gromov-Witten invariants, by
expanding the period map of the family $\chkX$ with respect to a
certain flat coordinate on the moduli space of complex structures of
mirror quintics.

As another example consider a two-torus $X$ (elliptic curve with
arbitrary choice of complex structure).  The $(1,1)$ Gromov-Witten
invariant of degree $d\geq 1$ with insertion the Poincar\'e dual class
of a point counts in this case the number of isogenies of degree $d$
to a fixed elliptic curve.  As such it satisfies
\[ \langle [\pt]^\PD \rangle_{1,1}^{X, d} = \sum_{k|d} k =
  \sigma_1(d), \]
and hence the generating series of these invariants (including the
$d=0$ case) is $-\frac{1}{24} E_2(Q)$ where $E_2$ denotes the
quasi-modular Eisenstein form of weight two.  The main result
of~(\cite{CalTu}) is that this equals the expansion in
$q = \exp(2\pi\ii\tau)$, around $q=0$, of the function of 
categorical enumerative $(1,1)$ invariants for the corresponding family
$\chkX$ of mirror elliptic curves.

\paragraph
Implicit in the above calculation for elliptic curves are the two
facts that
\begin{enumerate}
\item[(a)] $q$ is the flat coordinate, around the cusp, on the moduli
  space of elliptic curves; 
\item[(b)] the mirror map $\psi$ for elliptic curves identifies $q$
  with $Q$. 
\end{enumerate}

The main intuition behind Conjecture~\ref{conj:main} is a similar set
of assumptions, but for the flat coordinates around the hexagonal or
square points instead of around the cusp.  Below we will give precise
conjectural descriptions of the flat coordinates $q$ and $Q$
around the hexagonal point $\chkF\in M^\cx$ and its mirror
$F\in M^\Kah$.  The analysis for the square point is entirely similar.

\paragraph 
To understand these flat coordinates we need good descriptions of
$M^\Kah$ and $M^\cx$ around $F$ and $\chkF$.  We will review first the
classical situation (around the cusp) described in the work of
Polishchuk-Zaslow~\cite{PolZas}.

Polishchuk-Zaslow take the space $M^\Kah$ on a two-torus to be the
quotient of $\bbH$, with coordinate $\rho$, by $\rho \sim \rho+1$.
For each $\rho\in M^\Kah$ they construct a Fukaya category
$\cF^0(X^\rho)$ on the two-torus $X^\rho$ endowed with this structure.
The quotient above is precisely the same as the neighborhood of the
cusp on the moduli space $M^\cx$ of complex structures on a
two-torus\footnote{We ignore the stack structure of $M^\cx$, which
  only adds an extra $\bbZ/2\bbZ$ stabilizer.}.  For Polishchuk-Zaslow
the mirror map is simply the identity $\tau \leftrightarrow \rho$: the
complex elliptic curve $\chkX^\tau$ with modular parameter $\tau$
corresponds to the two-torus $X^\rho$ with complexified K\"ahler
structure $\rho = \tau$.

\paragraph
\label{subsec:mkah}
Even without explicitly constructing $M^\Kah$ as a moduli space of
geometric objects we could have understood its structure around the
large volume limit point through mirror symmetry.  Indeed, we could
have simply taken $M^\Kah$ to {\em be} the neighborhood of the large
complex structure limit point in $M^\cx$, a space we understand.  With
this point of view the mirror map is always the identity.

The same approach makes sense around the hexagonal point
$\chkF\in M^\cx$ and its mirror $F\in M^\Kah$.  The germ of $M^\cx$
around $\chkF$ is the quotient of $\bbH$ by
\[ \tau \sim \frac{\tau-1}{\tau}, \]
exhibiting the germ of $\bbH$ around $\tau_*$ as a triple cover of
$M^\cx$ branched over $\chkF$.  As above, we will {\em define} the
germ of $M^\Kah$ around $F$ to be the quotient of $\bbH$ (with
coordinate $\rho$) by $\rho\sim (\rho-1)/\rho$.  We think of
$\rho\in \bbH$ as giving an (ill-defined) ``complexified K\"ahler
class'' on the two torus, and write $X^\rho$ for this symplectic
geometry object.  The mirror map is, as before,
$\tau \leftrightarrow \rho$.

\paragraph
\label{subsec:justification}
In the A-model we conjecture that $Q = s(\rho)^3$ is a flat coordinate
on $M^\Kah$.  The justification for this comes from work of
Li-Shen-Zhou~\cite{LSZ}, where the authors suggest that the natural
way to interpret the generating series of FJRW invariants for two-tori
as a function of $\rho$ is via the map $s$ (with a different rescaling
from ours).  It would be natural to guess from their work that
$s(\rho)$ is the flat coordinate.  However, since $\rho$ is only defined up
to the equivalence $\rho \sim (\rho-1)/\rho$, the equality
\[ s\left (\frac{\rho-1}{\rho} \right )^3 = s(\rho)^3 \] implies that
$Q$ descends\footnote{This is not the only modification of $s(\tau)$
  that descends to a coordinate on $M^\Kah$, which in general will
  not be flat.  The same issue appears in the B-model.} to a
coordinate on $M^\Kah$, which we conjecture to be the flat
coordinate around $F$.

\paragraph
In the B-model Tu~\cite[Section 4]{Tu} argued that $h(t)/g(t)$ gives a
flat coordinate on the base $\bbA^1_t$ of the Hesse pencil of elliptic
curves,
\[ E_t:~\quad x^3+y^3+z^3+3txyz = 0. \]
Tu's work was motivated by a study of categories of graded matrix
factorizations, but via Orlov's correspondence~\cite{Orl} these are
equivalent to the derived categories of the above elliptic curves.

Again, $h(t)/g(t)$ does not give a coordinate on $M^\cx$ because
locally $\bbA^1_t$ is a branched triple cover of $M^\cx$ around
$\chkF$.  Its replacement $q = \left(h(t)/g(t)\right)^3$ does
descend to a coordinate on $M^\cx$ around $\chkF$, and we conjecture
it is flat.

\paragraph
By our construction of $M^\Kah$ the mirror map $\psi$ is the identity,
so the mirror of the complex curve $\chkX^\tau$ with modular parameter
$\tau$ is the symplectic object $X^\rho$ with
$\rho = \psi(\tau) = \tau$.  (Despite being equal we prefer to keep
$\rho$ and $\tau$ distinct since they represent different geometric
objects.)

Flat coordinates are unique up to multiplication by a scalar when the
moduli spaces $M^\Kah$ and $M^\cx$ are one-dimensional.  (The
rescaling factor $2\pi\Omega^2$ in~(\ref{subsec:rescale}) was chosen
so that this constant equals one.)  It follows that the flat
coordinates of $X^\rho$ and $\chkX^\tau$ are equal for $\rho = \tau$.

Consider a Hesse elliptic curve $E_t$ for some value of $t$.  It
can be written as $\chkX^\tau$ for some (non-unique) modular
parameter $\tau\in\bbH$.  The mirror of this curve is $X^\rho$ for
$\rho = \tau$.  (We think of $\rho\in M^\Kah$, so the ambiguity in
$\tau$ disappears.)  It follows that 
\[ \left(\frac{h(t)}{g(t)}\right)^3 = q(\chkX^\tau)  = Q(X^\rho) = s(\rho)^3, \]
or, using the fact that $s$ is invertible,
\[ s^{-1}\left( \frac{h(t)}{g(t)} \right ) \sim \rho \]
where $\sim$ is the equivalence relation used to define $M^\cx$
in~(\ref{subsec:mkah}).  Applying the $j$-function to both sides and
noting that it is $\sim$-invariant we get 
\[ j\left (s^{-1}\left( \frac{h(t)}{g(t)} \right ) \right ) = j(\rho)
  = j(E_t). \]
For the Hesse pencil the $j$-function can be computed
easily~\cite{ArtDol} and the result is
\[ j(E_t) = 27t^3\left (\frac{8-t^3}{1+t^3} \right )^3. \]
This is the statement of the conjecture.

\end{document}